\documentclass[12pt]{article}
\usepackage[francais]{babel}
\usepackage{amsmath,amsfonts,amssymb,amsthm}
\newcommand\A{\mathbb{A}}

\newcommand\C{\mathbb{C}}

\newcommand\R{\mathbb{R}}

\newcommand\Q{\mathbb{Q}}

\newcommand\dd{\mathfrak{d}}
\newcommand\Ss{\mathcal{S}}
\newcommand\bb{\mathcal{B}}

\newcommand\Oo{\mathcal{O}}

\newcommand\Bb{\mathcal{B}}

\newcommand\rg{\rightarrow}
\newcommand\lgr{\longrightarrow}

\newcommand\Res{\mathrm{Res}}

\def\adots{\mathinner{\mkern1mu\raise1pt\vbox{\kern7pt\hbox{.}}
\mkern2mu\raise4pt\hbox{.}
\mkern2mu\raise7pt\hbox{.}\mkern1mu}}

\begin{document}

\centerline{\LARGE Principe d'Heisenberg et fonctions positives}

\vskip4mm

\centerline{\large Jean Bourgain\footnote{Partially supported by NSF grant 0808042}\,,Ê Laurent Clozel\footnote{Membre de l'Institut Universitaire de France}  \ et Jean-Pierre Kahane}

\vskip 2cm

On d\'ecrit un probl\`eme naturel concernant les transform\'ees de Fourier, qui introduit des constantes $B_d$ d\'ependant de la dimension $d$.

Le probl\`eme est pos\'e, en dimension $1$, dans la premi\`ere section, o\`u on d\'ecrit une r\'eduction simple au cas des fonctions autoduales. Le premier r\'esultat est une minoration de $B_1$.

Puis on consid\`ere une classe tr\`es naturelle de fonctions qui donne simplement une majoration de $B_1$. Celle--ci  n'est pas optimale : on d\'ecrit des arguments simples qui permettent de l'am\'eliorer  (section~2).

Dans la section 3, ces calculs sont \'etendus aux dimensions arbitraires, donnant une minoration et une majoration simples des constantes.

Enfin, la section 4, arithm\'etique, relie ce probl\`eme \`a une question bien connue concernant les fonctions z\^{e}ta des corps de nombres.  Les arguments arithm\'etiques montrent  que la croissance lin\'eaire de $B_d$  en fonction de la dimension est naturelle au vu de propri\'et\'es connues de la ramification de ces corps.

\section{Position du probl\`eme et minoration de $B_1$}

Consid\'erons un couple de fonctions $(f,\widehat f)$ sur la droite r\'eelle ; c'est un couple de Fourier si
\[
\left\{
\begin{array}{ll}
 \widehat{f}(y) = \int f(x) e^{-2i\pi xy}dx\,, &f\in L^1 (\R)   \\
 \noalign{\vskip2mm}
 f(x) = \int\widehat{ f}(y) e^{2i\pi xy}dy\,, &\widehat f \in L^1 (\R)\,.   \\
\end{array}
\right.
\]

Ainsi $f$ et $\widehat f$ sont continues et tendent vers $0$ \`a l'infini. On s'int\'eressse aux couples de Fourier $(f,\widehat f)$ tels que

\begin{description}
\item[1)]  $f$ et $\widehat f$ sont r\'eelles et paires, non identiquement  nulles ;

\item[2)] $f(0) \le 0$ et $\widehat f (0)\le 0$ ;

\item[3)] $f(x)\ge 0$ pour $x\ge a_f$ et $\widehat f (y) \ge 0$ pour $y\ge a_{\widehat f}$.
\end{description}

Noter que la condition 2) et la non--nullit\'e de $f,\widehat f$ impliquent que $a_f$ et $a_{\widehat f}$ sont $>0$.

\vskip2mm

\noindent\textbf{Probl\`eme :} Quelle est la borne inf\'erieure du produit $a_fa_{\widehat f}$ pour les couples de Fourier $(f,\widehat f)$ v\'erifiant (1-3) ?

On d\'esignera cette borne inf\'erieure par $B_1\ge 0$ (noter que de tels couples existent \`a l'\'evidence). Nous allons montrer, ce qui n'est pas \'evident a priori, que  $B_1$ est strictement positif.

Jusqu'\`a la section 3, nous nous limiterons \`a la dimension 1. Pour un couple de Fourier v\'erifiant (1-3) posons
\[
\begin{array}{l}
  A(f) = \inf \{x >0 : f(]x,\infty[) \subset \R^+\}   \\
 A(\widehat f) = \inf \{y >0 : \widehat f (]y,\infty[) \subset \R^+\} \,.  \\
\end{array}
\]
Le produit $A(f)$ $A(\widehat f)$ est invariant par changement d'\'echelle, c'est--\`a--dire si on remplace $f(x)$, $\widehat f(y)$ par $f(x/\lambda)$, $\lambda \widehat f(y\lambda)$. Puisque
$$
B_1 =\inf \,A(f)\,A(\widehat f)
$$
pour tous les couples de Fourier v\'erifiant (1-3), on peut donc se limiter \`a ceux pour lesquels $A(f) = A(\widehat f)$. Alors $f+\widehat f\neq 0$ (consid\'erer ses valeurs en des points voisins de $A(f)$ et sup\'erieur \`a celui--ci), et
$$
A(f+\widehat f) \le A(f) = A(\widehat f)\,.
$$
Donc $B_1=\inf A^2(f+\widehat f)$. On voit donc que
$$
B_1 =A^2\, \qquad A=\inf A(f)
$$
la borne inf\'erieure \'etant prise sur l'ensemble des fonctions $f\in L^1 (\R)$, r\'eelles et paires, non identiquement nulles, \'egales \`a leur transform\'ee de Fourier, et telles que $f(0) \le 0$.

Posons
$$
\gamma(x)=e^{-\pi x^2}\,,
$$
de sorte que $\gamma=\widehat \gamma$. Si $f(0)<0$, $f-f(0)\gamma$ est non nulle et v\'erifie les m\^emes conditions que $f$, et
$$
A(f-f(0)\gamma) \le A(f)\,.
$$

Finalement,
$$
A= \inf A(f) \,, \leqno(1.1)
$$
\textbf{la borne inf\'erieure portant sur les fonctions $f\in L^1(\R)$ r\'eelles, paires, non identiquement nulles, et telles que $f=\widehat f$ et $f(0)=0$.}

\vskip2mm

Voici un r\'esultat important.

\vskip2mm

\textsc{Th\'eor\`eme} 1.  \textit{Soit $\lambda=-\inf\big({\sin x\over x}\big)=0,2172\cdots$}
$$
\begin{array}{lll}
\textit{Alors }  &A &\ge\displaystyle {1\over2(1+\lambda)}= 0,4107\cdots   \\
\noalign{\vskip2mm}
  \textit{donc }&B_1&\ge 0,1687\cdots   \\
 \end{array}
$$

\vskip 2mm

\textit{D\'emonstration}. Choisissons $f=\widehat f$, $f(0)=0$ et $\int_\R|f(x)|dx := \int_\R|f|=1$. Ecrivons simplement $A=A(f)$.  Posons $f=f^+-f^-$, $|f|=f^++f^-$. Comme $\int_\R f=\widehat f(0)=0$, on a $\int_\R f^+=\int_\R f^-=  \int_{-A}^{A} f^- = {1\over2}$. Donc $\int_{|x|\ge A}|f|=\int_{|x|\ge A}f^+ \le {1\over2}$, donc $\int_{|x|\le A}|f| \ge {1\over2}$. 
 Or $|f(x)|\le \int |\widehat f|=1$. Donc $2A\ge {1\over2}$, d'o\`u une premi\`ere minoration $A\ge {1\over4}$. On verra que cet argument s'\'etend aux dimensions sup\'erieures.
 
 En dimension 1, on peut le raffiner de la fa\c con suivante. Ecrivons, $f$ \'etant autoduale :
 $$
\begin{array}{c}
 \displaystyle  f(x)=\int f(y)\cos 2\pi yxdy=\int f(y)(\cos 2\pi yx-1)dy =   \\
 \noalign{\vskip2mm}
\displaystyle  = \int f^-(y) (1-\cos 2\pi yx)dy-\int f^+(y) (1-\cos 2\pi yx) dy\,.  
\end{array}
$$
Ceci implique, les deux int\'egrales \'etant positives :
$$
f^-(x) \le \int f^+ (y) (1-\cos 2 \pi yx)dy\,,
$$
d'o\`u
$$
{1\over4} = \int_0^A f^- \le \int_{-\infty}^{+\infty} f^+ (y) \big[A - {\sin 2\pi yA\over2\pi y}\big] dy
$$
et donc
$$
{1\over4} \le {A\over2}\sup_{u\in\R}\big(1-{\sin u\over u}\big) ={A\over2}(1+\lambda)
$$
d'o\`u le th\'eor\`eme.

\vskip2mm

Plus loin, nous aurons besoin de consid\'erer aussi des fonctions assez r\'eguli\`eres. Une classe naturelle est l'espace $\Ss$ de Schwartz. Il n'est point \'evident que la borne $A$ d\'efinie par (1.1), quand on impose de surcro\^{i}t \`a $f$ d'appartenir \`a $\Ss$, co\"{i}ncide avec celle d\'efinie pour $f$ parcourant $L^1$.

Notons $\Bb_1$ la constante $A^2$, o\`u $A$ est d\'efinie par (1.1) pour $f\in \Ss$. On va voir que $B_1$ et $\Bb_1$ diff\`erent assez peu. On a \'evidemment
$$
B_1 \le \Bb_1 \,. \leqno(1.2)
$$

Soit
$$
B_1^- = \inf(A^2\mid f(0) <0\,,\ f=\widehat f\ \textrm{paire}\neq0\,,\ f\in L^1)\,.
$$

Donc $B_1^-$ est d\'efinie par (1.1), o\`u l'on impose $f(0)<0$. On d\'efinit de m\^eme $\Bb_1^-$ en imposant de surcro\^{i}t  $f\in \Ss$.

A l'\'evidence :
$$
B_1^- \le \Bb_1^- \hskip 2,5cm\leqno(1.3)
$$
$$
\Bb_1 \le \Bb_1^-\,, \quad B_1\le B_1^-\,. \leqno(1.4)
$$

V\'erifions que $\Bb_1^- \le B_1^-$. Soit $f\in L^1$ v\'erifiant la condition (1.1) mais avec $f(0)<0$, et soit $a=A(f)$. Soit $\varphi =\psi * \psi$, $\psi$ \'etant $C^\infty$, paire, positive et de support compact tr\`es voisin de $0$, et $g=f*\varphi$. Alors $A(g) \le a+\varepsilon$ et $g(0)<0$. On a $\widehat g =\widehat f\widehat \psi^2$ ; en performant la m\^eme op\'eration sur $\widehat g$ on obtient une fonction $h\in \Ss$ telle que $h=\widehat h$, $h(0)<0$ et $A(h)\le a+\varepsilon$ ; on en d\'eduit que $\Bb_1^-\le B_1$ soit
$$
\Bb_1^- =B_1^-\,. \leqno(1.5)
$$

Noter que l'argument ne s'applique pas si $f(0)=0$. On va montrer
$$
B_1^- \le 2\, B_1\,; \leqno(1.6)
$$
d'apr\`es (1.4) et (1.6) on en d\'eduit
$$
B_1 \le \Bb_1\le 2\, B_1\,. \leqno(1.7)
$$

Soit $f$ v\'erifiant (1.1) et $a=A(f)$. Puisque $\widehat f(0)=\int f(x)dx=0$, $f$ prend des valeurs strictement n\'egatives sur $[-a,a]$. Soit $b>0$ tel que $f(b)<0$, et consid\'erons la distribution
$$
T= \delta _b +\delta _{-b} +2\delta _0\,.
$$

C'est une mesure positive, de type positif :
$$
\widehat T = 2\cos (2\pi by) +2\ge 0\,.
$$

On a
$$
(T*f)(0) =f(b)+f(-b) < 0\,.
$$

Puisque $b<a$, $g=T*f$ v\'erifie donc :
$$
g(0)<0\,,\quad g\ge 0\  \textrm{sur} \ [2a,\infty[\,.
 $$
 De plus $\widehat g =\widehat T \widehat f$ est $\ge 0$ sur $[a,\infty[$, et $ \widehat g(0)=0$.
 Par dilatation, on obtient alors une fonction $h$ telle que
 \[
\begin{array}{lll}
  h\ge 0&\textrm{sur}\ [a\sqrt{2}, \infty[\,,    &h(0) <0   \\
  \noalign{\vskip2mm}
 \widehat h\ge 0 &\textrm{sur}\ [a\sqrt{2}, \infty[\,,    &\widehat h(0) =0 \,.   \\ 
\end{array}
\]

Les fonctions $h$ et $\widehat h$ sont r\'eelles et paires. Alors $h+\widehat h$ v\'erifie les conditions relatives au calcul de $B_1^-$. Donc $B_1^-\le (a\sqrt{2})^2=2a$ ; variant $f$, on en d\'eduit enfin (1.6).

\section{Majoration de $B_1$}

Une premi\`ere id\'ee est d'associer \`a $f$ son d\'eveloppement d'Hermite
$$
f(x) \sim \sum_0^\infty a_n\ H_n(x)
$$
o\`u les $H_n$ sont des vecteurs propres de l'op\'erateur de Fourier $\mathcal {F}$, correspondant aux valeurs propres $i^n$. Ainsi $f=\widehat f$ s'exprime comme
$$
f(x) \sim\sum_0^\infty a_{4m} H_{4m}(x)\,.
$$

Les $H_n$ sont de la forme  $H_n(x)=e^{-\pi x^2}P_n(x)$ o\`u $P_n$ est un polyn\^ome de degr\'e $n$. Une combinaison lin\'eaire convenable de $H_0$ et $H_4$ (telle que $f(0)=0$) donne $\pi \, A^2\le 3$. Plus loin, les calculs semblent difficiles et nous n'avons pas poursuivi cette voie.

On peut aussi consid\'erer les fonctions
$$
g_a(x)= a\gamma (ax) +\gamma\bigg({x\over a}\bigg) - (1+a)\gamma(x)\,, \ a>1 \leqno(2.1)
$$
qui satisfont aux condition de (1.1). Alors toute expression de la forme
$$
\int_1^\infty g_a(x) d\tau(a) \leqno(2.2)
$$
o\`u $\tau$ est une mesure sur $]1,\infty[$ telle que l'int\'egrale converge pour tout $x$ et est $\ge 0$ \`a l'infini est une fonction candidate. (Il para\^{i}t difficile de d\'eterminer une propri\'et\'e simple et caract\'eristique de $\tau$ assurant que (2.2) est convergente et positive \`a l'infini).

Nous \'etudions d'abord  $A(g_a)$. Il est commode de poser $X=\pi x^2$, et $G_a(X)=g_a(x)$. Ainsi
$$
G_a(X) = a\, e^{-a^2X}+e^{-a^{-2}X} -(1+a) e^{-X}\,.
$$

De plus 
$$
H_a(X) = e^{X} G_a(X) =a\,e^{(1-a^2)X} +e^{(1-a^{-2})X}-1-a \leqno(2.3)
$$
est une fonction convexe, v\'erifiant
$$
H_a(0)=0\,, \quad H_a'(0) =-a^{-2}(a^2-1)(a^3-1) <0
$$
et tendant vers $+\infty$ avec $X$. Elle admet donc un unique z\'ero $X_a>0$, et
$$
A(g_a)=\sqrt{X_a\over\pi}\,.
$$
Il est naturel d'\'etudier la variation de $X_a$, et tout d'abord pour $a$ voisin de $1$. Posant $a=1+h$ ,  $h>0$, il vient pour $X$ fix\'e
$$
H_a(X)=(1+h)(e^{-X(2h+h^2)}-1) +e^{X(2h-3h^2+3h^3-4h^4)X}-1
$$
modulo $O(h^5)$. Ceci s'\'ecrit $P_1h+P_2h^2+ P_3h^3 +P_4h^5+O(h^5)$, o\`u les polyn\^omes $P_i$ sont :
\[
\begin{array}{l}
 P_1 =0   \\
  P_2=2X(2X-3)   \\
 P_3=-X(2X-3) \\
 \displaystyle P_4= -5X+15X^2 -{28\over3} X^3 +{4\over3} X^4\,. 
\end{array}
\]

L'expression de $P_2$ montre que pour $h$ assez petit $H_a(X)>0$ si $X>{3\over2}$ et $H_a(X) <0$ si $X<{3\over2}$. Par cons\'equent,
$$
\lim_{a\rg 1^+} X_a ={3\over2}\,. \leqno(2.4)
$$
Ce qui fournit une borne explicite
$$
A\le \sqrt{3\over2\pi}\,. \leqno(2.5)
$$

Mais cette borne simple ne peut \^etre la vraie valeur de $A$. Pour $X={3\over2}$, $P_2$ et $P_3$ s'annulent, et
$$
P_4\bigg({3\over2}\bigg)={3\over2}\,.
$$
Pour $h$ petit et non nul, on a donc $X_a<{3\over2}$.

Si $a\rg +\infty$, $X_a\rg +\infty$ ; en fait, un calcul simple montre que
$$
X_a=\log a+O(1)\qquad (a\rg +\infty)\,.
$$

Nous n'avons pas d\'etermin\'e la valeur minimale de $X_a$, mais il est facile de l'estimer, de fa\c con semi--heuristique. La valeur $a=\sqrt{2}$ donne, en posant $q=e^{{1\over2}X_a}$ :
$$
q^3 - (1+\sqrt{2}) q^2 +\sqrt{2}=0\, ;
$$
si $q\neq 1$, l'\'equation quadratique
$$
q^2-\sqrt{2}q -\sqrt{2}=0
$$
donne $q={\sqrt{2}\over2} (1+\sqrt{1+2\sqrt{2}})$,
$$
X_a=2\,\log q=1,4749\ldots <{3\over2}\quad (a=\sqrt{2})\,.
$$
La valeur $a=2$ donne, pour $q=e^{{3\over4}X}$ :
$$
q^4 -2 {q^4-1\over q-1}=0\,.
$$
La solution $q>1$ est $q=2,9744\ldots$, d'o\`u
$$
X_a=1,4534\ldots \qquad \qquad (a=2)\,.
$$
Il est vraisemblable que c'est \`a peu pr\`es la valeur optimale accessible par cette m\'ethode. En effet si l'on r\'esout $H_a(X)=0$, $H_a$ donn\'ee par (2.3), et que l'on suppose $a\ge 2$, le premier terme est n\'egligeable. Donc $X_a$ est \`a peu pr\`es
$$
\log (1+a)\over1-a^{-2}\,.
$$
L'extremum de cette expression est atteint pour $a(1-a)=2$ log$(1+a)$, qui donne
$$
a=2,08137\ldots
$$

Dans tous les cas, la valeur minimale $A(g_a)$ ainsi obtenue n'est pas la valeur (1.1) cherch\'ee. Consid\'erons en effet $a_0$ tel que $X_0=X_{a_0}$ soit minimal, et $H_0=H_{a_0}$, positive sur $[X_0,\infty[$.

Soit $a$ (par exemple, voisin de 1) tel que $X_a>X_0$. Sur $[X_a,\infty[$, $H_a$ est $\ge 0$ et son ordre de croissance pour $X\rg +\infty$, en $e^{(1-a^{-2})X}$, est plus petit que celui de $H_{a_0}$ si $a<a_0$. Il existe donc $T>0$ tel que $H_{a_0}-TH_a$ soit $\ge 0$ sur $[X_a,\infty[$. Mais cette fonction est $>0$ sur $[X_0,X_a[$, donc sur un voisinage de $X_0$, donc pour $X\ge X'$ avec $X'<X_0$.

Le m\^eme argument s'applique en prenant tout $a_0$ tel que $X_0<{3\over2}$. Pour $a_0=2$, on peut d\'eterminer la correction optimale (qui correspond \`a $a$ tr\`es voisin de 1), donnant une fonction $\ge 0$ sur $[X'',\infty[$, avec
$$
\begin{array}{ll}
X'' &= 1,25\ldots\\
A&\le 0,63\ldots
\end{array}
\leqno(2.6)
$$
Nous n'avons fait qu'un calcul tr\`es approch\'e. Enon\c cons n\'eanmoins le r\'esultat, \`a comparer au th\'eor\`eme~1.

\vskip2mm

\textsc{Th\'eor\`eme} 2. \textit{On a $A \le 0,64$ et $B_1\le 0,41$.}

\section{Dimensions sup\'erieures}

Nous nous pla\c cons  dans $\R^d$ euclidien, produit scalaire
$$
x\cdot y = \sum_1^d x_iy_i\,, \qquad \|x\| =(x\cdot x)^{1/2}\,,
$$
la transform\'ee de Fourier \'etant donn\'ee par 
$$
\widehat f (y) =\int f(x) e^{-2i\pi x\cdot y}dx \leqno(3.1)
$$
o\`u $dx =dx_1\ldots dx_d$ est la mesure de Lebesgue ; alors
$$
f(x) = \int\widehat{ f}(y)e^{2i\pi x\cdot y}dy\,. \leqno(3.2)
$$

On suppose $f,\widehat f$ continues et int\'egrables. Plus g\'en\'eralement, si $E$ est un espace euclidien de dimension $d$, si la mesure invariante $dx$ sur $E$ est choisie de sorte que la mesure du cube engendr\'e par une base orthonormale soit \'egale \`a 1, et si $x\cdot y$ d\'esigne le produit scalaire, la transform\'ee de Fourier (3.1) a pour r\'eciproque (3.2).

On consid\`ere les couples de Fourier $(f,\widehat f)$ v\'erifiant (3.3)

\begin{description}
\item[1)] $f$, $\widehat f$ r\'eelles et paires, non identiquement nulles
\item[2)] $f(0)\le 0$ et $\widehat f(0)\le 0$
\item[3)] $f(x)\ge 0$ pour $\|x\|\ge a_f$, $\widehat f(y)\ge0$ pour $\|y\|\ge a_{\widehat f}$. 
\end{description}

On d\'efinit, comme dans le \S 1, $A(f)$ et $A(\widehat f)$ :
$$
A(f) = \inf\{r>0 : f(x)\ge 0Ê \ \ \textrm{si}\ \ \|x\| >r\}\,,
$$
et
$$
B_d =\inf A(f)A(\widehat f)
$$
pour les couples v\'erifiant 1), 2), 3). Soit $f^\#(x)$ l'int\'egrale (invariante) de $f$ sur la sph\`ere de rayon $\|x\| : \widehat f^\#=(\widehat f)^\#$, et $f^\#$ et $\widehat f^\#$ ne sont pas nulles ; sinon $f$ et $\widehat f$ seraient \`a support compact d'apr\`es   3). Puisque $A(f^\#)\le A(f)$ et $A(\widehat f^\#) \le A(\widehat f)$, on peut se limiter aux couples de fonctions radiales. Puisque
$$
(f(x/\lambda))^\wedge = \lambda^d \widehat f(\lambda y)\qquad \qquad (\lambda>0)\,,
$$
l'argument du \S 1 s'applique alors et l'on voit que
$$
B_d=A^2\,, \qquad A=\inf A(f)\,, \leqno(3.4)
$$
\textbf{la borne inf\'erieure \'etant prise sur l'ensemble des fonctions $f\in L^1(\R^d)$, radiales, non identiquement nulles, et telles que $f=\widehat f$ et $f(0)=0$.}

On a, comme dans le \S 1, ajout\'e si n\'ecessaire un multiple de la fonction, radiale et autoduale
$$
\gamma(x) = e^{-\pi\|x\|^2}\,.
$$

\vskip2mm

\textsc{Th\'eor\`eme} 3. \textit{On a $B_d\ge {1\over\pi}\big({1\over2}\Gamma \big({d\over2}+1\big)\big)^{2/d}> {d\over2\pi e}$}

\vskip2mm

\textit{D\'emonstration.} Elle est calqu\'ee sur le cas $d=1$, en rempla\c cant l'intervalle $(-A(f),A(f))$ par la boule de centre $O$ et de rayon $A(f)$, dont le volume $(\ge {1\over2})$ est ${1\over\Gamma({d\over2}+1)}(A(f))^d \pi^{d/2}$.

\vskip2mm

Posant $X=\pi\|x\|^2$, l'argument du \S 2 nous am\`ene \`a consid\'erer les fonctions naturelles
$$
g_a(x)=G_a(X)\qquad (x\in \R^d)
$$
o\`u
$$
G_a(X) =a^d e^{-X a^2}+e^{-Xa^{-2}}-(1+a^d) e^{-X}\,,
$$
et enfin
$$
H_a(X)=a^d e^{(1-a^2)X}+e^{(1-a^{-2})X}- (1+a^d)\,, \ a>1\,.
$$
Il est commode de poser $a^2=1+k$, $d=2c$, d'o\`u
$$
H_a(X) =(1+k)^c \,e^{-kX} +e^{(1-(1+k)^{-1})X} -1-(1+k)^c.    
$$
La d\'eriv\'ee \`a l'origine en $X$ est
$$
{k\over1+k}\bigg(1-(1+k)^{c+1}\bigg) <0\,;
$$
l'argument de convexit\'e du \S 2 montre que $H_a$ a un unique z\'ero positif $X_a$. Comme auparavant, nous calculons un d\'eveloppement en $k$  limit\'e \`a l'ordre 4 de $H_a(X)$ . Il vient
$$
H_a(X)=P_1k+P_2k^2+P_3k^3+P_4k^4+O(k^5)\,,
$$
$$
\begin{array}{ll}
  P_1   &=0   \\
  P_2   &=X(X-c-1)   \\
  P_3   &=\displaystyle{1\over2}(c-2)X(X-c-1)\\
  \noalign{\vskip2mm}
  P_4 &=\displaystyle{1\over12}X \{
  X^3-(2c+6)X^2+(3c(c-1)+18)X-\\
  &\hfill -(2c(c-1)(c-2)+12) \}.
\end{array}
$$

Comme en dimension 1, on voit que $P_2$ et $P_3$ s'annulent pour
$$
X=X(d) := {d\over2}+1 \leqno(3.4)
$$
De plus $P_2$ est $>0$ pour $X>X(d)$\
,   $<0$ pour $X<X(d)$. Faisant tendre $k$ vers $0$ on en d\'eduit
$$
\lim_{a\rg 1} X_a = {d\over2}+1
$$

Pour comprendre la position de $X_a$ par rapport \`a $X(d)$ quand $a\rg 1$, calculons $Q_4$(X(d))
, o\`u $P_4={X\over12}Q_4$. Le calcul donne
$$
Q_4(c+1) =-c^2+1\,.
$$

Pour $d>2$, ce terme est donc $<0$, donc $H_a(X(d))<0$ pour $a$ proche de 1, ce qui montre que 
$$
X_a>{d\over2}+1\qquad (a>1\,,\ \textrm{assez proche de }\ 1)\,.
$$

Il est donc possible que la valeur (3.4) soit optimale. Pour $d=1$, ce n'est pas le cas, comme on l'a vu au~\S 2.

Pour $d=2$
, $Q_4(c+1)=0$, donc nous devons calculer, \`a l'ordre 5 au moins, le d\'eveloppement limit\'e de
$$
H_a(2)=(1+k)e^{-2k}+e^{2(1-{1\over1+k})} -2-k \,. \leqno(3.5)
$$
Le d\'eveloppement de Taylor en $0$ de
\[
\begin{array}{ll}
f(z)   &=e^{2(1-{1\over1+z})} = e^{2{z\over1+z}} :  \\
  f(z)  &=\displaystyle \sum_0^\infty q_n\ z^n\,,   \\ 
\end{array}
\]
se calcule par le th\'eor\`eme des r\'esidus. Posant
$$
w={z\over1+z}\,,\ z={w\over1-w}\,, \ dz={dw\over(w-1)^2}\,,
$$
il vient, les int\'egrales \'etant prises sur un petit contour autour de $0$ :
\[
\begin{array}{ll}
  q_n  &= \Res_{z=0} \displaystyle{f(z)\over z^{n+1}}= {1\over2i\pi} \oint e^{2z\over1+z} {dz\over z^{n+1}}   \\
  \noalign{\vskip2mm}
  &= \displaystyle{1\over2i\pi} \oint e^{2w}{(1-w)^{n+1}\over w^{n+1}} {dw\over(1-w)^2}   \\
  \noalign{\vskip2mm}
  &= \Res_{w=0}\displaystyle{(1-w)^{n-1}\over w^{n+1}} e^{2w}\,.   
\end{array}
\]

En particulier, $p_5$ est la somme de
$$
{2^4\over4!} - {2^5\over5!} \leqno(3.6)
$$
venant du premier terme de (3.5), et du terme en $w^5$ de $e^{2w}(1-w)^4$, \'egal \`a
$$
{2^5\over5!} - 4\cdot {2^4\over4!} +6\cdot {2^3\over3!} - 4\cdot {2^2\over2!} +2\,. \leqno(3.7)
$$

On trouve que $p_5=0$.

De m\^eme, $p_6$ est la somme de
$$
-{2^5\over5!} + {2^6\over6!} \leqno(3.8)
$$
et de
$$
{2^6\over6!} -5\cdot {2^5\over5!} +10\cdot {2^4\over4!} -10\cdot {2^3\over3!} +5\cdot {2^2\over2!}- 2\,, \leqno(3.9)
$$
d'o\`u
$$
p_6=-{4\over45} <0\,.
$$

Pour $a$ tr\`es voisin de 1, on a donc $H_a(2)<0$ et $X_a > X(2)=2$. L\`a encore, il est possible que la borne donn\'ee par (3.4) soit optimale.

Pour conclure ce paragraphe, noter que l'on a obtenu pour tout $d\ge 2$ la borne sup\'erieure
$$
B_d\le \Bb_d\le {d+2\over2\pi}\leqno(3.10)
$$
o\`u $\Bb_d$ est d\'efini, comme dans le \S 1, par les fonctions de l'espace $\Ss(\R^d)$. Par ailleurs les consid\'erations de la fin du \S 1, relatives aux bornes pour $L^1$ et pour $\Ss$, s'appliquent. Dans la d\'emonstration de l'in\'egalit\'e (1.6), on doit consid\'erer $T=\delta _b + \delta _{-b}+2\delta _0$, o\`u $\|b\|<a=A(f)$ et $f(b)<0$ ; $\widehat T=2\cos(2\pi\, b\cdot y)+2$ est une onde plane positive. Le reste de l'argument est identique, en rempla\c cant $h+\widehat h$ par la moyenne sph\'erique de $h+\widehat h$ si on veut s'en tenir aux fonctions radiales. En conclusion, \`a comparer au Th\'eor\`eme~3 :

\vskip2mm

\textsc{Th\'eor\`eme} 4. \textit{On a}
$$
B_d\le \Bb_d \le {d+2\over2\pi}\,,\quad B_d\ge {1\over2}\Bb_d \,. \leqno(3.11)
$$

\section{Un argument arithm\'etique}

Soit $F$ un corps de nombres de degr\'e $d$ sur $\Q$. On d\'esigne par $v$ les places (finies ou archim\'ediennes) de $F$, et par $F_v$ la compl\'etion correspondante ; pour $v$ finie $\Oo_v\subset F_v$ est l'anneau des entiers et $\Oo_v^\times$ son groupe des unit\'es ; $q_v$ est le cardinal du corps r\'esiduel. Soit
$$
\A_F = \prod_v{}' F_v
$$
(produit restreint) l'anneau des ad\`eles de $F$, et $\A_F^\times =I_F$ le groupe des id\`eles. Soit $|\ | : x \in I_F \mapsto \prod_v |x|_v$ la norme d'id\`ele,
$$
\begin{array}{cll}
  &I_F^1   &=\{x\in I_F : |x|=1\}  \\
  \noalign{\vskip2mm}
\textrm{et} &I_F^+   &=\{x\in I_F : |x|\ge1\}\,.  \\
\end{array}
$$

On consid\`ere la mesure invariante $dx=\prod dx_v$ sur $\A_F$, $dx_v$ \'etant  une mesure  de Haar sur $F_v$. En les places finies, $dx_v$ est la mesure autoduale de Tate \cite{tate} ; en une place r\'eelle ; $dx$ est la mesure de Lebesgue ; en une place complexe, dont l'on note $z=x+iy$ la variable, $dz=2dxdy$. En une place r\'eelle, la transform\'ee de Fourier $\widehat f(y)$ d'une fonction $f$ est d\'efinie comme dans le reste de cet article.

Si $z=x+iy$ est le param\`etre en une place complexe, et $w=\xi+i\eta$, Tate d\'efinit la transform\'ee $\widehat f(w)$ d'une fonction $f(z)$ par
\[
\begin{array}{lll}
  &\widehat f(w)   &=\displaystyle \int f(z) e^{-2i\pi Tr(zw)}dz   \\
  \noalign{\vskip2mm}
 \textrm{o\`u} &Tr(zw)   &= 2\textrm{Re}(zw) = 2(x\xi-y\eta)\,.   \\
\end{array}
\]

Pour des fonctions radiales, donc paires en chacune des variables, ceci revient \`a consid\'erer la transform\'ee de Fourier d\'efinie, comme dans le \S 3, par le produit scalaire $z\cdot w= 2(x\xi+y\eta)$. La mesure autoduale $dz$ de Tate est la mesure normalis\'ee consid\'er\'ee au d\'ebut du \S 3 pour un espace euclidien abstrait.

Soit $f$ la fonction de l'espace de Schwartz de $\A_F$ donn\'ee par 
$$
f(x) = \prod_{v|\infty} f_v(x_v) \prod_{v\textrm{ finie}} f_v^0(x_v) \leqno(4.1)
$$
o\`u $f_v^0$ est la fonction caract\'eristique de $\Oo_v$ et o\`u, pour $v$ archim\'edienne, $f_v$ est pour l'instant une fonction arbitraire de l'espace de Schwartz. La fonction z\^eta de Tate associ\'ee est d\'efinie pour $\textrm{Re}(s)>1$ par
$$
Z(f,s) = \int_{I_F} f(x) |x|^s\ d^\times x\,,
$$
o\`u $d^\times x$ est le produit des $d^\times x_v$, $d^\times x_v= {dx_v\over|x_v|}$  (multipli\'e par $(1-q_v^{-1})^{-1}$ aux places finies).

Plut\^ot que les fonctions d\'ecompos\'ees de (4.1), nous consid\'erons, sur $\R^d$, des fonctions de la forme $g_a(x)$ (\S 3) o\`u $\R^d$ est consid\'er\'e comme un espace euclidien par
$$
\|x_\infty\|^2 = \sum_{v\textrm{ r\'eelle}} |x_v|^2 + \sum_{v\textrm{ complexe}} 2\|z_v\|^2\,,
$$
$\|z\|$ \'etant la valeur absolue usuelle d'un nombre complexe. (On notera $|z|=\|z\|^2$ la valeur absolue normalis\'ee comme dans la th\'eorie de Tate). Plus g\'en\'eralement,
$$
f(x) =f_\infty(x_\infty) \prod_{v\textrm{ finie}} f_v^0(x_v) \leqno(4.2)
$$
o\`u $f_\infty(x_\infty)\in \Ss(\R^d)$. Les conditions impos\'ees par Tate (i.e., $(z
_1),\, ( z_2)\,,(z_3)$ in \cite{tate}, \S 4.4) sont v\'erifi\'ees par de telles fonctions. Par exemple, $(z_3)$ prescrit que  l'int\'egrale
$$
\int_{F_\infty} f_\infty (x_\infty) \prod_{v\mid \infty} |x_v|_v^{\sigma{-1}} dx\,,
$$
o\`u $F_\infty=\prod\limits_{v\mid\infty}F_v$, soit absolument convergente pour $\sigma>1$. C'est vrai en fait pour $\sigma>0$ et tout $f_\infty\in \Ss(F_\infty)$. La m\^eme condition est donc v\'erifi\'ee pour~$\widehat f$.

Dans le cas o\`u $f_\infty=\prod f_v^0$, avec 
\[
\begin{array}{ll}
  f_v^0(x) = e^{-\pi x^2}   &\textrm{(variable r\'eelle)}   \\
  f_v^0(z)= e^{-2\pi\|z\|^2} &\textrm{(variable complexe)\,,}   \\
\end{array}
\]
$Z(f,s)$ est la fonction z\^eta $\zeta_F(s)$, multipli\'ee par ses facteurs archim\'ediens usuels (produit de fonctions $\Gamma$). Ecrivons, d'apr\`es Tate,
$$
\begin{array}{lrl}
  Z(f,s)   &= &\displaystyle\int_{I_F^+} f(x) |x|^s d^\times x  +\int_{I_F^+} \widehat f(x)|x|^{1-s}d^\times x \\
  \noalign{\vskip2mm}
  &   &\displaystyle+\kappa{\widehat f(0)\over s-1} - \kappa{f(0)\over s}   \\   
\end{array}
\leqno(4.3)
$$
o\`u, avec les notations usuelles (\cite{tate}, Th\'eor\`eme~4.3.2),
$$
\kappa= {2^{r_1}(2\pi)^{r_2}hR\over\sqrt{D_F}\ w}
$$
est le r\'esidu en $s=1$ de $\zeta_F(s)$. En particulier, $D_F$ est le discriminant absolu de $F$, et $d=r_1+2r_2$, $r_1$ \'etant le nombre de places r\'eelles et $r_2$ le nombre de places complexes. Les deux int\'egrales figurant dans (4.3) sont alors absolument convergentes pour tout $s\in \C$.

\vskip2mm

\noindent\textsc{Lemme} 1.--- 
\textit{Soit $s$ un z\'ero de $\zeta_F(s)$ tel que $\textrm{Re}(s)>0$. Alors $Z(f,s)$ s'annule en $s$ pour tout choix de $f_\infty \in \Ss(F_\infty)$.}

\vskip2mm

En effet on peut \'ecrire, d'abord pour $\textrm{Re}s>1$,
$$
Z(f,s) = Z(f_\infty,s)\zeta_F(s)
\,.$$

\vskip2mm

 Puisque $Z(f,s)$, ainsi que $\zeta_F(s)$ et $Z(f_\infty,s)$ sont holomorphes pour $s\neq1$, $\textrm{Re}(s)>0$, le Lemme s'en d\'eduit.

\vskip2mm

Pour toute place finie $v$, $\widehat f_v^0$ est \'egale \`a $|\mathfrak{d}_v|^{-1/2}\textrm{char}(\mathfrak{d}_v^{-1})$. Ici $\dd_v\subset F_v$ est la diff\'erente, $\dd_v^{-1}$ son inverse, $\textrm{char}(\dd_v^{-1})$ la fonction caract\'eristique, et $|\dd_v|$ est la norme d'id\'eal (une puissance positive de $q_v$). Rappelons que
$$
\prod_{v\textrm{ finie}} |\dd_v| = |D_F|\,.
$$

Consid\'erons alors la premi\`ere int\'egrale de (4.3) : 
$$
\int_{I_F^+} f(x) \ |x|^s d^\times x\,.
\leqno(4.4)
$$

Si $f(x) \neq 0$ en $x=(x_\infty,x_f)$, la description de $f_f=\prod\limits_{v\textrm{ finie}}f_v$ montre que $|x_{f}|\le 1$ ; puisque $|x_\infty x_f| \ge 1$,
$$
 |x_\infty| = \prod_{v\mid\infty}|x_v| \ge 1\,.
\leqno(4.5)
$$

Dans la deuxi\`eme int\'egrale, en remarquant que $|x_v|\le |\dd_v|$ si $x_v\in \dd_v^{-1}$, on a de m\^eme $|x_{f}| \le \prod\limits_v|\dd_v| = |D_F|$ d'o\`u
$$
|x_\infty| \ge D_F^{-1}\,.
\leqno(4.6)
$$

\vskip2mm

\noindent\textsc{Lemme} 2.--- \textit{Supposons qu'il existe un couple de Fourier $(f,\widehat f)$ sur $F_\infty=\R^d$ tel que $f(x_\infty)\ge 0$ si $\mid x_\infty|\ge 1$, $f$ prend des valeurs strictement positives au voisinage de $1$ dans l'ensemble $|x_\infty|\ge 1$, $\widehat f(y_\infty)\ge 0$ si $|y_\infty| \ge D_F^{-1}$ et $f(0)=\widehat f(0)=0$. Alors $\zeta_F(s)\neq 0$ pour tout $s$ dans l'intervalle $]0,1[$.}

\vskip2mm

En effet (4.3) est alors r\'eduit \`a ses termes int\'egraux ; $|x|^s$ est strictement positif dans le domaine d'int\'egration, et l'int\'egrale (4.4) est strictement positive d'apr\`es la propri\'et\'e impos\'ee \`a $f$. Donc $Z(f,s)>0$ et $\zeta_F(s)\neq0$ d'apr\`es le Lemme~1.

\vskip2mm

Soit $x=(x_v)\in F_\infty$. La norme euclidienne compatible avec la transform\'ee de Fourier de Tate est
$$
\|x\|^2= \sum_{v\textrm{ r\'eelle}} |x_v|^2 + 2 \sum_{v \textrm{ complexe}} \|x_v\|^2\,.
$$
Puisque
$$
|x|^2 = \prod_{v \textrm{ r\'eelle}} |x_v|^2 \prod_{v\textrm{ complexe}} \|x_v\|^4\,,
$$
l'in\'egalit\'e arithm\'etico--g\'eom\'etrique donne
$$
|x|^{2/d}\le {1\over d} \|x\|^2\,.
$$

Posant $r=\|x\|$, $\rho=\|y\|$ $(y\in F_\infty)$ on voit que 
\[
\begin{array}{lll}
|x|\ge 1  &\Longrightarrow   & r\ge \sqrt{d}   \\
|y| \ge |D_F|^{-1}  &\Longrightarrow   &\rho\ge |D_F|^{-1/d} \sqrt{d}\,.   \\
\end{array}
\]

\vskip2mm

\noindent\textsc{Proposition} 1.--- \textit{Supposons qu'il existe un corps de nombres $F$ de degr\'e $d$ et de discriminant $D$ tel que $\zeta_F(s)$ a un z\'ero dans $]0,1[$. Alors}
$$
\Bb_d\ge d\ |D|^{-1/d}\,.
$$

\vskip2mm

R\'eciproquement, bien s\^ur, $\zeta_F$ n'a pas de z\'ero r\'eel si
$$
d\ |D|^{-1/d} > \Bb_d\,.
$$
La d\'emonstration est maintenant \'evidente. Supposons pour exemple que $d\ |D|^{-1/d}>\bb_d$. On peut trouver $f$, $\widehat f$ radiales, comme dans le \S 3, et $\ge0$ pour $r\ge \sqrt{d}$ et $\rho\ge |D|^{-1/d}\sqrt{d}$. On peut supposer aussi que $f$ prend des valeurs strictement positives sur l'ensemble des $x$ tels que $\sqrt{d}\le \|x\| \le \sqrt{d}+\varepsilon$. Les conditions du Lemme~2 sont alors r\'eunies puisque $\|1\|=\sqrt{d}$.

\vskip2mm

Il est difficile de trouver des corps $F$ v\'erifiant l'hypoth\`ese de la Proposition. Cependant, la d\'ecomposition, pour $F$ galoisien sur $E$, de $\zeta_F(s)$ en fonctions $L$ d'Artin pour $E$ a permis \`a Armitage d'exhiber un tel z\'ero (en $s={1\over2}$ bien s\^ur, conform\'ement \`a l'hypoth\`ese de Riemann).

Plus pr\'ecis\'ement, Armitage consid\`ere une extension explicite $F$ de $E=\Q(\sqrt{3(1+i)})$, de degr\'e 12 sur $E$ et donc de degr\'e 48 sur $\Q$, construite par Serre \cite{serre}, et montre que $\zeta_F\big({1\over2}\big)=0$. (\cite{armi}, \S 5). 

Par cons\'equent, la th\'eorie des nombres impliquait a priori la version faible suivante du Th\'eor\`eme~3 :

\vskip2mm

\noindent\textsc{Proposition} 2.--- \textit{Si $d$ est multiple de $48$, $\bb_d$ est strictement positif.}

\vskip2mm

Pour $d=48$, ceci r\'esulte de l'existence de $F$. Supposons que $d=48c$. Il existe une extension cyclotomique $L$ de $\Q$ de degr\'e $c$ et lin\'eairement disjointe de $F$. Alors $LF$ est de degr\'e $d$ sur $\Q$, et $\zeta_F$ divise $\zeta_{LF}$ puisque $LF/F$ est ab\'elienne, et que $\zeta_{LF}$ se factorise donc en produit de fonctions $L$ de Dirichlet relatives \`a $F$. D'o\`u le r\'esultat.

\vskip2mm

On peut se demander si la Proposition 1 implique une restriction sur les discriminants des corps tels que $\zeta_F$ ait un z\'ero r\'eel. Dans ce cas, on~a
$$
|D|^{1/d} \ge {d\over\bb_d} \,. \leqno(4.7)
$$

Mais, inconditionnellement d'apr\`es le Th\'eor\`eme 3,
$$
{d\over\bb_d} < 2\pi e = 17,079\cdots
$$

Or les minorations d'Odlyzko \cite{odly} donnent en g\'en\'eral
$$
  |D|^{1/d} \ge 22,2 (1+0(d))
$$
pour $d\lgr \infty$. Par cons\'equent (4.7) est automatiquement v\'erifi\'e, au moins pour $d$ assez grand.

La proposition 2 ne conduit donc pas \`a une minoration int\'eressante de $\bb_d$. Il est frappant de remarquer n\'eanmoins que, pour certains degr\'es au moins, la Th\'eorie des nombres impliquait la croissance lin\'eaire en $d$ donn\'ee par le Th\'eor\`eme~3.
 Soit en effet $p$ un nombre premier. D'apr\`es les th\'eor\`emes de Golod--Shafarevi\v{c} et Brumer, il existe une suite de corps
$$
E_p^1 \subset E_p^2 \subset \cdots E_p^n \subset \cdots
$$
o\`u $E_p^1$, de degr\'e $p(p-1)$ sur $\Q$, est une extension de degr\'e $p$ de $\Q(\zeta_p)$ et o\`u $E_p^{n+1}/E_p^n$ est ab\'elienne, non ramifi\'ee de degr\'e $p$. Voir \cite{roqu}, Cor.~7 ; on a adjoint $\zeta_p$ pour obtenir $E_p^1$ par deux extensions successives, ab\'eliennes, \`a partir de~$\Q$.

Consid\'erons la suite d'extensions $F_i=F\, E_p^i$ de $F$, $F_i/F_{i+1}$ \'etant ab\'elienne, de degr\'e 1 ou $p$. On peut en extraire une sous--suite minimale strictement croissante, d'o\`u
$$
F_0 = F\, E_p^{n_0} \subset F_1 \subset \cdots \subset F_m = F\, E_p^{n_m} \subset \cdots\,,
$$
chaque extension ab\'elienne de degr\'e $p$. Vu l'absence de ramification relative, une formule classique donne l'expression des discriminants absolus :
$$
D(F_m) = D(F_0)^{p^m}:= D^{p^m}\,.
\leqno(4.8)
$$

Les extensions successives \`a partir de $F$ \'etant ab\'eliennes, $\zeta_F$ divise $\zeta_{F_m}$ pour tout $m$. La Proposition~1 donne alors pour $d=d_0 p^m$, $d_0=[F_0:\Q]$ :
$$
\bb_d \ge C\, d\,,\quad C=|D|^{-1_/d_0}\,.
\leqno(4.9)
$$

Pour de telles suites de degr\'es, (3.10) et (4.8) montrent donc que la croissance de $\bb_d$ --- et donc de $B_d\ge {1\over2}\bb_d$ --- est lin\'eaire en $d$. Si $p$ ne divise pas $D_F$, $F$ et $\Q(\zeta_p)$ sont lin\'eairement disjoints et l'on peut choisir $E_p^1$ lin\'eairement disjoint de $F$. Alors $F_0=F\, E_p^1$ et l'in\'egalit\'e (4.8) est valide pour $d=48(p-1)p^n$, $n\ge 1$. Bien s\^{u}r, le terme en $(p-1)$ n'est pas n\'ecessaire si l'on est pr\^{e}t \`a utiliser la conjecture d'Artin ou la conjecture de divisibilit\'e de Dedekind.

\eject

\vskip2mm

\hfill\begin{minipage}{6,5cm}
Jean Bourgain

School of Mathematics

Institute for Advanced Study,

Princeton, N.J. 08540 

Etats Unis

\vskip2mm

Laurent Clozel, Jean--Pierre Kahane

Laboratoire de Math\'ematique

Universit\'e Paris--Sud, B\^at. 425

91405 Orsay Cedex

France

\end{minipage}

\end{document}